\documentclass[12pt,twoside]{article}
\usepackage{amsmath,amssymb}    
\usepackage{url}                
\usepackage{graphicx}           
\usepackage{amsthm}
\theoremstyle{definition}
\newtheorem{theorem}{Theorem}[section]
\newtheorem{lemma}[theorem]{Lemma}
\newtheorem{proposition}[theorem]{Proposition}
\newtheorem{corollary}[theorem]{Corollary}

\newtheorem{definition}[theorem]{Definition}

\newtheorem{remark}[theorem]{Remark}

\newcommand{\C}{\mathbb{C}}

\date{}
\newcommand{\h}{{\mathfrak{h}}}

\begin{document}

\pagestyle{myheadings}
\title{A uniform proof of the Macdonald-Mehta-Opdam
identity for finite Coxeter groups}
\author{Pavel Etingof}
\maketitle

\section{Introduction}

In this note we give a new proof of the 
Macdonald-Mehta-Opdam integral identity for finite Coxeter 
groups. This identity was conjectured by Macdonald and proved by Opdam 
in \cite{O1,O2} using the theory of multivariable Bessel functions,
but in non-crystallographic cases the proof relied on a 
computer calculation by F. Garvan. 
Our proof is somewhat more elementary (in particular,
it does not use multivariable Bessel functions), and 
uniform (does not refer to the classification of 
finite Coxeter groups).  

{\bf Acknowledgements.} I am very grateful to Ivan Cherednik,
whose explanations regarding shifting the contour of integration 
led me to the main idea of this proof. I would also like to thank Misha Feigin 
and Charles Dunkl for reading the preliminary version of this note and making useful comments. 
The work of the author was  partially supported by the NSF grant
DMS-0504847.
\section{Preliminaries}

\subsection{Coxeter groups}

Let $W$ be a finite  
Coxeter group of rank $r$ with reflection representation $\h_{\Bbb R}$
equipped with a Euclidean $W$-invariant inner product $(,)$.
\footnote{As a basic reference on finite Coxeter groups, we use the book \cite{Hu}.} 
Denote by $\h$ the complexification of $\h_{\Bbb R}$. 
The reflection hyperplanes subdivide $\h_{\Bbb R}$ into $|W|$ 
chambers; let us pick one of them to be the 
dominant chamber and call its interior $D$. 
For each reflection hyperplane, pick the perpendicular vector 
$\alpha\in \h_{\Bbb R}$ with $(\alpha,\alpha)=2$ 
which has positive inner products with elements of $D$, and 
call it the positive root corresponding to this hyperplane.
The walls of $D$ are then defined by the equations $(\alpha_i,v)=0$, 
where $\alpha_i$ are simple roots. 
Denote by $S$ the set of positive roots, 
and for $\alpha\in S$ denote by $s_\alpha$ the corresponding reflection. 
We will denote the set of reflections also by $S$. Let 
$$
\Delta(x)=\prod_{\alpha\in S}(\alpha,x)
$$ 
be the corresponding discriminant 
polynomial. Let $d_i,i=1,...,r$, be the degrees of the
generators of the algebra $\Bbb C[\h]^W$. 
Note that $|W|=\prod_i d_i$. 

\subsection{Cherednik algebras}

For $k\in \C$, let $H_k=H_k(W)$ be the corresponding rational Cherednik 
algebra (see e.g. \cite{E}). Namely, $H_k$ 
is the quotient of $\C[W]\ltimes T(\h\oplus \h)$ 
(with the two generating copies of $\h$ spanned by $x_a,y_a$, $a \in \h$), 
by the defining relations
$$
[x_a,x_b]=[y_a,y_b]=0, [y_a,x_b]=(a,b)+
k\sum_{\alpha\in S}(\alpha,a)(\alpha,b)s_\alpha.
$$

Let $a_i$ be an orthonormal basis of $\h$. 
Consider the element
$$
\bold h=\sum_i x_{a_i}y_{a_i}+\frac{r}{2}
+k\sum_{\alpha\in S}s.
$$
It satisfies $[\bold h,x_a]=x_a, [\bold h,y_a]=-y_a$.

Let $M_k=H_k\otimes_{\Bbb CW\ltimes \Bbb C[y_{a_i}]}\Bbb C$, where 
$y_{a_i}$ act in $\Bbb C$ by $0$ and $w\in W$ by $1$. 
Then we have a natural vector space isomorphism $M_k\cong \Bbb C[\h]$.
For this reason $M_k$ is called {\em the polynomial representation} of $H_k$. 
The elements $y_{a_i}$ act in this representation by Dunkl operators (see \cite{E}). 

\begin{proposition}\label{contra} There exists a unique
$W$-invariant symmetric bilinear form $\beta_k$ on $M_k$
such that $\beta_k(1,1)=1$, which
satisfies the contravariance condition
$$
\beta_k(y_av,v')=\beta_k(v,x_av'),\ v,v'\in M_k, a\in \h.
$$
Polynomials of different degree are orthogonal under $\beta_k$. 
Moreover, the kernel of $\beta_k$ is the maximal proper submodule 
of $M_k$, so $M_k$ is reducible iff $\beta_k$ is degenerate.  
\end{proposition}

\begin{proof}
The proof is standard, see e.g. \cite{E}.
Namely, let $M_k^*$ be the dual space of $M_k$ 
with the dual action of $H_k$ twisted by the antiautomorphism
of $H_k$ given by $x_a\to y_a$, $y_a\to x_a$, and $w\to w^{-1}$, 
$w\in W$. Then a symmetric $W$-invariant 
bilinear form $\beta: M_k\times M_k\to \Bbb C$
is the same thing as an $H_k$-homomorphism 
$\hat\beta: M_k\to M_k^*$. 
Since this homomorphism commutes with ${\bold h}$, it must land in the 
graded dual space $M_k^\dagger\subset M_k^*$ and preserve the grading. 
But such a homomorphism clearly exists and is unique up to scaling, as it is 
determined by $\hat\beta(1)$. 
This implies the existence and uniqueness of $\beta_k$, 
and the fact that polynomials of different degrees 
are orthogonal under $\beta_k$. 

Now, it is clear from the definition that the 
kernel of $\beta_k$ is a submodule in $M_k$, 
so it remains to show that the module $M_k/{\rm Ker}\beta_k$ is 
irreducible. For this, let $L_k$ be the irreducible quotient of 
$M_k$; then we have a natural surjective homomorphism $M_k\to L_k^\dagger$
(defined up to scaling), which must factor through $L_k$. 
Thus we have a diagram 
$$
M_k\to L_k\cong L_k^\dagger\to M_k^\dagger,
$$ 
which implies that $\hat \beta_k$ factors through $L_k$, 
i.e. $M_k/{\rm Ker}\beta_k=L_k$, as desired. 
\end{proof}

\section{The main theorem}

The goal of this note is to give a uniform and self-contained proof 
of the following theorem. 

\begin{theorem}\label{maint}
(i) (The Macdonald-Mehta integral) For ${\rm Re}(k)\ge 0$, one has 
$$
(2\pi)^{-r/2}\int_{\h_{\Bbb R}}e^{-(x,x)/2}|\Delta(x)|^{2k}dx=
\prod_{i=1}^r \frac{\Gamma(1+kd_i)}{\Gamma(1+k)}.
$$

(ii) Let $b(k):=\beta_k(\Delta,\Delta)$. Then
$$
b(k)=|W|\prod_{i=1}^r\prod_{m=1}^{d_i-1}(kd_i+m).
$$
\end{theorem}

For Weyl groups, this theorem was proved by E. Opdam \cite{O1}.
The non-crystallographic cases were done by Opdam in \cite{O2} 
using a direct computation in the rank 2 case 
(reducing (i) to the beta integral), 
and a computer calculation by F. Garvan for $H_3$ and $H_4$. 

In the next subsection, we give a uniform proof of Theorem \ref{maint}.
We emphasize that many parts of this proof are borrowed from  
Opdam's previous proof of this theorem. 

\section{Proof of the main theorem}

\begin{proposition}\label{c1}
The function $b$ is a polynomial of degree at most $|S|$, and 
the roots of $b$ are negative rational numbers. 
\end{proposition} 

\begin{proof}
Since $\Delta$ has degree $|S|$,  
it follows from the definition of $b$ 
that it is a polynomial of degree $\le |S|$. 

Suppose that $b(k)=0$ for some $k\in \Bbb C$. Then $\beta_k(\Delta,P)=0$ for any polynomial $P$.
Indeed, if ${\rm deg}(P)\ne |S|$, this follows from Proposition \ref{contra},
while if $P$ has degree $|S|$, this follows from the fact that $\Delta$ is the 
unique (up to scaling) polynomial of degree $|S|$ that is antisymmetric under $W$.  
 
Thus, $M_k$ is reducible and hence has a singular vector, i.e. 
a nonzero homogeneous polynomial $f$ of positive degree $d$ 
living in an irreducible representation $\tau$ of $W$ killed by $y_a$. 
Applying the element ${\bold h}$ to $f$, we get 
$$
k=-\frac{d}{m_\tau},
$$
where $m_\tau$ is the eigenvalue of the operator $T:=\sum_{\alpha\in S}(1-s_\alpha)$ 
on $\tau$. But it is clear (by computing the trace of $T$) 
that $m_\tau\ge 0$. This implies that any root of $b$ is negative rational. 
\end{proof}

Denote the Macdonald-Mehta integral by $F(k)$. 

\begin{proposition}\label{l1}
One has 
$$
F(k+1)=b(k)F(k).
$$
\end{proposition}

\begin{proof}
Let $\bold f=\frac{1}{2}\sum y_{a_i}^2$.
Introduce the {\em Gaussian inner product} on $M_k$ as follows:

\begin{definition}
The Gaussian inner product
$\gamma_k$ on $M_k$
is given by the formula
$$
\gamma_k(v,v')=\beta_k(\exp(\bold f)v,\exp(\bold
f)v').
$$
\end{definition}

This makes sense because the operator $\bold f$ is locally nilpotent on
$M_k$. 

Note that $\Delta$ is a nonzero $W$-antisymmetric polynomial of the 
smallest possible degree, so $(\sum y_{a_i}^2)\Delta=0$ and hence 
\begin{equation}\label{bk}
\gamma_k(\Delta,\Delta)=\beta_k(\Delta,\Delta)=b(k).
\end{equation} 

\begin{proposition}\label{x-inv}
Up to scaling, $\gamma_k$
is the unique $W$-invariant symmetric bilinear form on $M_k$ satisfying
the condition
$$
\gamma_k((x_a-y_a)v,v')=\gamma_k(v,y_av'),\ a\in
\h.
$$
\end{proposition}

\begin{proof} We have
$$
\gamma_k((x_a-y_a)v,v')=
\beta_k(\exp(\bold f)(x_a-y_a)v,\exp(\bold
f)v')=
$$
$$
\beta_k(x_a\exp(\bold f)v,\exp(\bold
f)v')=
\beta_k(\exp(\bold f)v,y_a\exp(\bold
f)v')=
$$
$$
\beta_k(\exp(\bold f)v,\exp(\bold
f)y_av')=
\gamma_k(v,y_av').
$$

 Let us now show uniqueness.
If $\gamma$ is any $W$-invariant
symmetric bilinear form satisfying the condition of the Proposition, then
let $\beta(v,v')=\gamma(\exp(-\bold f)v,\exp(-\bold f)v')$.
Then $\beta$ is contravariant, so by Proposition \ref{contra},
it's a multiple of $\beta_k$, hence $\gamma$ is a multiple of
$\gamma_k$.
\end{proof}

Now we will need the following known result (see \cite{Du2}, Theorem 3.10).

\begin{proposition}\label{inte}
For ${\rm Re}(k)\ge 0$ we have
\begin{equation}\label{intformu}
\gamma_k(f,g)=F(k)^{-1}\int_{\h_{\Bbb R}}f(x)g(x)d\mu_c(x)
\end{equation}
where
$$
d\mu_c(x):=e^{-(x,x)/2}|\Delta(x)|^{2k}dx.
$$
\end{proposition}

\begin{proof}
It follows from Proposition \ref{x-inv} that $\gamma_k$
is uniquely, up to scaling,
determined by the condition that it is $W$-invariant,
and $y_a^\dagger=x_a-y_a$. These properties
are easy to check for the right hand side of 
(\ref{intformu}), using the fact that
the action of $y_a$ is given by Dunkl operators.
\end{proof}

Now we can complete the proof of Proposition \ref{l1}.
By Proposition \ref{inte}, we have  
$$
F(k+1)=F(k)\gamma_k(\Delta,\Delta),
$$
so by (\ref{bk}) we have 
$$
F(k+1)=F(k)b(k).
$$
\end{proof}

Let 
$$
b(k)=b_0\prod (k+k_i)^{n_i}.
$$
We know that $k_i>0$, and also $b_0>0$
(because the inner product $\beta_0$ on real polynomials 
is positive definite).  

\begin{corollary}\label{c2}
We have 
$$
F(k)=b_0^k\prod_i \left(\frac{\Gamma(k+k_i)}{\Gamma(k_i)}\right)^{n_i}.
$$
\end{corollary}

\begin{proof} Denote the right hand side by $F_*(k)$ 
and let $\phi(k)=F(k)/F_*(k)$. Clearly, $\phi(1)=1$.
Proposition \ref{l1} implies that $\phi(k)$ is a 1-periodic 
positive function on $[0,\infty)$. Also by the Cauchy-Schwarz inequality, 
$$
F(k)F(k')\ge F((k+k')/2)^2,
$$
so $\log F(k)$ is convex for $k\ge 0$. 
This implies that $\phi=1$, since $(\log F_*(k))''\to 0$ as $k\to +\infty$. 
\end{proof} 

In particular, we see from Corollary \ref{c2} and the multiplication formulas for the $\Gamma$ function 
that part (ii) of the main theorem implies part (i).

It remains to establish (ii). 

\begin{proposition}\label{degree}
The polynomial $b$ has degree exactly $|S|$. 
\end{proposition}

\begin{proof}
By Proposition \ref{c1}, $b$ is a polynomial of degree at most $|S|$.
To see that the degree is precisely $|S|$, let us 
make the change of variable 
$y=k^{1/2}x$ in the Macdonald-Mehta integral 
and use the steepest descent method. 
We find that the leading term of the asymptotics 
of $\log F(k)$ as $k\to +\infty$ is 
$|S|k\log k$. This together with 
the Stirling formula and Corollary \ref{c2} 
implies the statement.    
\end{proof}

\begin{proposition}\label{l2}
The function 
$$
G(k):=F(k)\prod_{j=1}^r \frac{1-e^{2\pi ikd_j}}{1-e^{2\pi ik}}
$$
analytically continues to an entire function of $k$.
\end{proposition}

\begin{proof}
Let $\xi\in D$ be an element. Consider the real hyperplane 
$C_t=it\xi+\Bbb \h_{\Bbb R}$, $t>0$. Then $C_t$ does not intersect 
reflection hyperplanes, so we have a continuous branch 
of $\Delta(x)^{2k}$ on $C_t$ which tends to the positive branch in $D$ as 
$t\to 0$. Then, it is easy to see that 
for any $w\in W$, the limit of this branch in 
the chamber $w(D)$ will be  
$e^{2\pi ikl(w)}|\Delta(x)|^{2k}$. 
Therefore, by letting $t=0$, we get 
$$
(2\pi)^{-r/2}\int_{C_t}e^{-(x,x)/2}\Delta(x)^{2k}dx=
\frac{1}{|W|}F(k)(\sum_{w\in W}e^{2\pi ikl(w)})
$$
(as this integral does not depend on $t$). 
But it is well known that 
$$
\sum_{w\in W}e^{2\pi ikl(w)}=
\prod_{j=1}^r \frac{1-e^{2\pi ikd_j}}{1-e^{2\pi ik}},
$$
(\cite{Hu}, p.73), so
$$ 
(2\pi)^{-r/2}|W|\int_{C_t}e^{-(x,x)/2}\Delta(x)^{2k}dx=G(k).
$$
Since  $\int_{C_t}e^{-(x,x)/2}\Delta(x)^{2k}dx$ is clearly an entire function,
the statement is proved. 

\end{proof}

\begin{corollary}\label{c3}
For every $k_0\in [-1,0]$ 
the total multiplicity of all the roots of $b$
of the form $k_0-p$, $p\in \Bbb Z_+$, equals
the number of ways to represent $k_0$ in the form $-m/d_i$, 
$m=1,...,d_i-1$. In other words,  
the roots of $b$ are $k_{i,m}=-m/d_i-p_{i,m}$, $1\le m\le d_i-1$, 
where $p_{i,m}\in \Bbb Z_+$. 
\end{corollary}

\begin{proof}
We have 
$$
G(k-p)=\frac{F(k)}{b(k-1)...b(k-p)}
\prod_{j=1}^r \frac{1-e^{2\pi ikd_j}}{1-e^{2\pi ik}},
$$
Now plug in $k=1+k_0$ and large positive integer $p$. 
Since by Proposition \ref{l2} the left hand side is regular, so must be the right hand side, 
which implies the claimed upper bound for the total multiplicity, 
as $F(1+k_0)>0$. The fact that the bound is actually attained follows from 
the fact that the polynomial $b$ has degree exactly $|S|$ (Proposition \ref{degree}), and the fact that 
all roots of $b$ are negative rational (Proposition \ref{c1}). 
\end{proof}

It remains to show that in fact in Corollary \ref{c3}, $p_{i,m}=0$ for all $i,m$;
this would imply (ii) and hence (i). 

\begin{proposition}\label{k2}
Identity (i) of the main theorem is satisfied in $\Bbb C[k]/k^2$. 
\end{proposition}

\begin{proof}
Indeed, we clearly have $F(0)=1$. Next, a rank $1$ computation gives
$F'(0)=-\gamma|S|$, where $\gamma$ is the Euler constant, while the  
derivative of the right hand side of (i) at zero equals to 
$$
-\gamma\sum_{i=1}^r (d_i-1).
$$
But it is well known that 
$$
\sum_{i=1}^r (d_i-1)
=|S|,
$$ 
(\cite{Hu}, p.62), 
which implies the result.
\end{proof}

\begin{remark}
In fact, Proposition \ref{k2} allows one to make 
Opdam's original proof of the main theorem given in \cite{O2} 
classification independent and computer-free. Indeed, 
the arguments of \cite{O2} imply that (i) holds 
up to a factor of the form $c^k$, where $c>0$, 
and Proposition \ref{k2} implies that $c=1$. 
\end{remark}

\begin{proposition}\label{l5} Identity (i) of the main theorem is satisfied in $\Bbb C[k]/k^3$. 
\end{proposition}

Note that Proposition \ref{l5} immediately implies (ii), and hence the whole theorem. Indeed, 
it yields that 
$$
(\log F)''(0)=\sum_{i=1}^r\sum_{m=1}^{d_i-1}(\log \Gamma)''(m/d_i),
$$
so by Corollary \ref{c3} 
$$
\sum_{i=1}^r\sum_{m=1}^{d_i-1}(\log \Gamma)''(m/d_i+p_{i,m})=
\sum_{i=1}^r\sum_{m=1}^{d_i-1}(\log \Gamma)''(m/d_i),
$$
which implies that $p_{i,m}=0$ since $(\log \Gamma)''$ is strictly decreasing on $[0,\infty)$. 

\begin{proof} (of Proposition \ref{l5})
We will need the following result 
about finite Coxeter groups.
Let $\psi(W)=3|S|^2-\sum_{i=1}^r (d_i^2-1)$. 

\begin{lemma}\label{rk2}
One has 
\begin{equation}\label{rk2e}
\psi(W)=\sum_{G\in {\rm Par}_2(W)}\psi(G),
\end{equation}
where ${\rm Par}_2(W)$ is the set of parabolic subgroups of $W$ of rank 2. 
\end{lemma}
 
\begin{proof}
Let 
$$
Q(q)=|W|\prod_{i=1}^r\frac{1-q}{1-q^{d_i}}.
$$  
It follows from Chevalley's theorem that
$$
Q(q)=(1-q)^r\sum_{w\in W}\det(1-qw|_\h)^{-1}.
$$
Let us subtract the terms for $w=1$ and $w\in S$ from both sides of this equation, 
divide both sides by $(q-1)^2$, and set $q=1$ (cf. \cite{Hu}, p.62, formula (21)). Let $W_2$ be the set of elements of $W$
that can be written as a product of two different reflections. Then by a straightforward computation 
we get
$$
\frac{1}{24}\psi(W)=\sum_{w\in W_2}\frac{1}{r-{\rm Tr}_\h(w)}. 
$$
In particular, this is true for rank 2 groups. 
The result follows, as any element $w\in W_2$ belongs 
to a unique parabolic subgroup $G_w$ of rank $2$
(namely, the stabilizer of a generic point $\h^w$, 
\cite{Hu}, p.22). 
\end{proof}

Now we are ready to prove the proposition.  
By Proposition \ref{k2}, it suffices to show the coincidence of the second derivatives of (i) 
at $k=0$. The second derivative of the right hand side of (i) at zero 
is equal to 
$$
\frac{\pi^2}{6}\sum_{i=1}^r (d_i^2-1)+\gamma^2 |S|^2.
$$
On the other hand, we have 
$$
F''(0)=(2\pi)^{-r/2}\sum_{\alpha,\beta\in S}\int_\h e^{-(x,x)/2}\log \alpha^2(x)\log \beta^2(x)dx.
$$
Thus, from a rank 1 computation we see that our job is to 
establish the equality 
$$
(2\pi)^{-r/2}\sum_{\alpha\ne \beta\in S}\int_\h e^{-(x,x)/2}\log \alpha^2(x)\log \frac{\beta^2(x)}{\alpha^2(x)}dx
$$
$$
=\frac{\pi^2}{6}(\sum_{i=1}^r (d_i^2-1)-3|S|^2)=-\frac{\pi^2}{6}\psi(W). 
$$
Since this equality holds in rank 2 (as in this case (i) reduces to the beta integral),  
in general it reduces to equation (\ref{rk2e}) (as for any $\alpha\ne \beta\in S$,
$s_\alpha$ and $s_\beta$ are contained in 
a unique parabolic subgroup of $W$ of rank 2). 
The proposition is proved. 
\end{proof}

\end{document}